\documentclass[11pt]{amsart}
\usepackage{mathrsfs}
\usepackage{cases}
\usepackage{amssymb}
\usepackage{amsfonts}
\usepackage{amssymb,amsfonts}

\usepackage{amsmath,amssymb,amsthm,hyperref}
\allowdisplaybreaks

\begin{document}
\theoremstyle{plain}
\newtheorem{thm}{Theorem}[section]
\newtheorem{theorem}[thm]{Theorem}
\newtheorem{lemma}[thm]{Lemma}
\newtheorem{corollary}[thm]{Corollary}
\newtheorem{corollary*}[thm]{Corollary*}
\newtheorem{proposition}[thm]{Proposition}
\newtheorem{proposition*}[thm]{Proposition*}
\newtheorem{conjecture}[thm]{Conjecture}
\theoremstyle{definition}
\newtheorem{construction}{Construction}
\newtheorem{notations}[thm]{Notations}
\newtheorem{question}[thm]{Question}
\newtheorem{problem}[thm]{Problem}
\newtheorem{remark}[thm]{Remark}
\newtheorem{remarks}[thm]{Remarks}
\newtheorem{definition}[thm]{Definition}
\newtheorem{claim}[thm]{Claim}
\newtheorem{assumption}[thm]{Assumption}
\newtheorem{assumptions}[thm]{Assumptions}
\newtheorem{properties}[thm]{Properties}
\newtheorem{example}[thm]{Example}
\newtheorem{comments}[thm]{Comments}
\newtheorem{blank}[thm]{}
\newtheorem{observation}[thm]{Observation}
\newtheorem{defn-thm}[thm]{Definition-Theorem}

\newcommand{\sM}{{\mathcal M}}


\title{Integrality of the LMOV invariants for framed unknot}
       \author{Wei Luo, Shengmao Zhu}
        \address{Center of Mathematical Sciences, Zhejiang University, Hangzhou, Zhejiang 310027, China}
        \email{luowei@cms.zju.edu.cn, szhu@zju.edu.cn}

\begin{abstract}
The Labastida-Marin\~o-Ooguri-Vafa (LMOV) invariants are the open
string BPS invariants which are expected to be integers based on the
string duality conjecture from M-theory. Several explicit formulae
of LMOV invariants for framed unknot have been obtained in the
literature. In this paper, we present a unified method to deal with
the integrality of such explicit formulae. Furthermore, we also
prove the integrality of certain LMOV invariants for framed unknot
in higher genera.
\end{abstract}
\maketitle

\section{Introduction}
Topological string amplitude is the generating function of
Gromov-Witten invariants which are usually rational numbers
according to their definitions \cite{HKKPTVVZ}. In 1998, Gopakumar
and Vafa \cite{GV0,GV1} found that topological string amplitude is
also the generating function of a series of integer-valued
invariants related to BPS counting in M-theory. Later, Ooguri and
Vafa \cite{OV} extended the above result to open string case, the
corresponding integer-valued invariants are named as OV invariants.
Furthermore, the OV invariants are further refined by Labasitida,
Mari\~no and Vafa in \cite{LMV}, then the resulted invariants are
called LMOV invariants \cite{LP}. An expanded physicist's
reconsideration of the GV and LMOV can be found in \cite{DW}. We
refer to \cite{Zhu3} for a brief review of the applications of these
integrality structures of topological strings in mathematics.

The open string LMOV invariants have been studied in many papers,
such as
\cite{LM1,LM2,LMV,LP,CLPZ,GKS,KS1,LZ,MMMS,MMMRSS,KRSS1,KRSS2}. Based
on the large $N$ duality of topological string and Chern-Simons
theory \cite{W,GV2,OV}, the open string LMOV invariants can be
approached by investigating the colored HOMFLYPT invariants of the
dual knots. For a knot $\mathcal{K}$, we use the notation
$n_{\mu,g,Q}(\mathcal{K})$ to denote the LMOV invariants of genus
$g$ with a boundary type $\mu$ which is a partition of a positive
integer, where $Q$ is a parameter describing the dependence of the
relative homology class of the dual Calabi-Yau geometry of the knot
$\mathcal{K}$. We refer to Section \ref{section2} for more detailed
definition of $n_{\mu,g,Q}(\mathcal{K})$.

In particular, when the knot $\mathcal{K}$ is a framed unknot
$U_\tau$ with framing $\tau\in \mathbb{Z}$, we have different ways
to compute its LMOV invariants according to string dualities which
have been proved in this situation \cite{LLZ,Zhou,EO}.
 Then one can obtain several
explicit formulae \cite{GKS,LZ} for the genus zero LMOV invariants
$n_{\mu,0,Q}(U_\tau)$ of framed unknot $U_\tau$. It turns out that
these explicit formulae are certain combinations of the M\"obius
function and binomial numbers. Based on the integrality conjecture
for LMOV invariants, these formulae will give integers. However,
such an argument is not so obvious, a rigorous proof is required.

In this paper, we present a straightforward way to prove the
integrality of these formulae.  We use the notation $n_{m,l}(\tau)$
to denote the LMOV invariants $n_{(m),0,l-\frac{m}{2}}(U_\tau)$ of
the framed unknot $U_\tau$ of genus $0$, where $m\geq 1$ and $l\geq
0$ are two integers. We have the following explicit formula
\cite{MV,LZ} for $n_{m,l}(\tau)$.

For $b\geq 0$ and $a\in \mathbb{Z}$, we introduce the notation
$\binom{a}{b}$ which is defined as follows
 \makeatletter
\let\@@@alph\@alph
\def\@alph#1{\ifcase#1\or \or $'$\or $''$\fi}\makeatother
\begin{subnumcases}
 {\binom{a}{b}=} 1, & $b=0$, \nonumber  \\\nonumber \binom{a}{b},
&$b\geq 1$ and $a\ge 0$,\\\nonumber (-1)^{b}\binom{-a+b-1}{b},
&$b\geq 1$ and $a<0$.
\end{subnumcases}
\makeatletter\let\@alph\@@@alph\makeatother
We define
\begin{align*}
c_{m,l}(\tau)=-\frac{(-1)^{m\tau+m+l}}{m^2}\binom{m}{l}\binom{m\tau+l-1}{m-1},
\end{align*}
then
\begin{align} \label{integralformula}
n_{m,l}(\tau)=\sum_{d|m,d|l}\frac{\mu(d)}{d^2}c_{\frac{m}{d},\frac{l}{d}}(\tau),
\end{align}
where $\mu(d)$ denotes the M\"{o}bius functions.

In Section \ref{section3}, we prove that
\begin{theorem} \label{Thm1}
For any $\tau\in \mathbb{Z}$, $m\geq 1, l\geq 0$, we have
$n_{m,l}(\tau) \in \mathbb{Z}$.
\end{theorem}

\begin{remark}
In fact, such form of the formula (\ref{integralformula}) is very
general. For example, if we take some special values of $l$ or
$\tau$, it will give the formulae in \cite{GKS} (cf. the formulae
(1.4) and (1.5) in \cite{GKS}):
\begin{equation} \label{GKS1}
    b_{K_p,r}^{-}=-\frac{1}{r^2}\sum_{d\mid r}\mu(\frac{r}{d}) \binom{3d-1}{d-1}, \qquad
    b_{K_p,r}^{+}=\frac{1}{r^2}\sum_{d\mid r}\mu(\frac{r}{d}) \binom{(2|p|+1)d-1}{d-1}
\end{equation}
for $p\leq -1$ and
\begin{equation} \label{GKS2}
    b_{K_p,r}^{-}=-\frac{1}{r^2}\sum_{d\mid r}\mu(\frac{r}{d}) (-1)^{d+1}\binom{2d-1}{d-1}, \qquad
    b_{K_p,r}^{+}=\frac{1}{r^2}\sum_{d\mid r}\mu(\frac{r}{d}) (-1)^{d} \binom{(2p+2)d-1}{d-1}
\end{equation}
for $p\geq 2$. The above formulae (\ref{GKS1}) and (\ref{GKS2}) are
referred as the extremal BPS invariants of twist knots in
\cite{GKS}. Therefore, Theorem \ref{Thm1} implies the integrality of
formulae (\ref{GKS1}) and (\ref{GKS2}) immediately. Moreover, the
integrality of another special case of the formula
(\ref{integralformula}) was also proved in \cite{Zhu2}.
\end{remark}

Then, denoted by $n_{(m_1,m_2)}(\tau)$ the LMOV invariants
$n_{(m_1,m_2),0,\frac{m_1+m_2}{2}}(U_\tau)$ of the framed unknot
$U_\tau$ with  $\mu=(m_1,m_2)$, $g=0$ and $Q=\frac{m_1+m_2}{2}$,
where $m_1\geq m_2\geq 1$, we obtain the following formula
\begin{align} \label{integralformula2}
    n_{(m_1,m_2)}(\tau)& =  \frac{1}{m_1+m_2}\sum_{d\mid m_1,d\mid m_2} \mu(d) (-1)^{(m_1+m_2)(\tau+1)/d}  \\
                     & \cdot
                     \binom{(m_1\tau+m_1)/d-1}{m_1/d}\binom{(m_2\tau+m_2)/d}{m_2/d}\nonumber.
\end{align}
From this expression, we know that
$n_{(m_1,m_2)}(\tau)=n_{(m_2,m_1)}(\tau)$. With the similar method,
in Section \ref{section4}, we prove that
\begin{theorem} \label{Thm2}
For any $\tau\in \mathbb{Z}$ and $m_1,m_2\geq 1$, then
$n_{(m_1,m_2)}(\tau)\in \mathbb{Z}$.
\end{theorem}

Next, let $n_{m,g,Q}(\tau)$ be the LMOV invariants
$n_{(m),g,Q}(U_\tau)$ of higher genus $g$ with boundary condition
$\mu=(m)$. We define the following generating function
\begin{align*}
g_m(q,a)=\sum_{g\geq 0}\sum_{Q\in
\mathbb{Z}/2}n_{m,g,Q}(\tau)z^{2g-2}a^Q
\end{align*}
where $z=q^{\frac{1}{2}}-q^{-\frac{1}{2}}$.

Let
\begin{align*}
\mathcal{Z}_m(q,a)=(-1)^{m\tau}\sum_{|\nu|=m}\frac{1}{\mathfrak{z}_\nu}\frac{\{m\nu\tau\}}{\{m\}\{m\tau\}}\frac{\{\nu\}_a}{\{\nu\}}
\end{align*}
where $\mathfrak{z}_{\nu}=|Aut(\nu)|\prod_{i=1}^{l(\nu)}\nu_i$ and
$\{m\}$ denotes the quantum integer, see Section \ref{section2} for
 introduction of the above notations.

By the definition of LMOV invariants in Section \ref{section2}, we
obtain the following expression
\begin{align*}
g_m(q,a)=\sum_{d|m}\mu(d)\mathcal{Z}_{m/d}(q^d,a^d).
\end{align*}

In Section \ref{section5}, we prove that
\begin{theorem} \label{Thm3}
For any $m\geq 1$, we have
$
g_m(q,a)\in z^{-2}\mathbb{Z}[z^2,a^{\pm \frac{1}{2}}],
$
where $z=q^{\frac{1}{2}}-q^{-\frac{1}{2}}$.
\end{theorem}
Therefore, Theorem \ref{Thm3} implies that $ n_{m,g,Q}(\tau)\in
\mathbb{Z} $ and moreover $n_{m,g,Q}(\tau)$ vanishes for large $g$
and $Q$.

\textbf{Acknowledgements.}
 The authors are grateful to the referees for careful reading of the paper and
 valuable comments and suggestions
 which greatly improved the presentation of the content.

\section{LMOV invariants} \label{section2}
\subsection{Basic notations}
We first introduce some basic notations. A partition $\lambda$ is a finite sequence of positive integers $%
(\lambda_1,\lambda_2,..)$ such that $\lambda_1\geq
\lambda_2\geq\cdots$. The length of $\lambda$ is the total number of
parts in $\lambda$ and denoted by
$l(\lambda)$. The weight of $\lambda$ is defined by $|\lambda|=%
\sum_{i=1}^{l(\lambda)}\lambda_i$. The automorphism group of
$\lambda$, denoted by Aut($\lambda$), contains all the permutations
that
permute parts of $\lambda$ by keeping it as a partition. Obviously, Aut($%
\lambda$) has the order $|\text{Aut}(\lambda)|=\prod_{i=1}^{l(\lambda)}m_i(%
\lambda)! $ where $m_i(\lambda)$ denotes the number of times that
$i$ occurs in $\lambda$. Define
$\mathfrak{z}_{\lambda}=|\text{Aut}(\lambda)|\prod_{i=1}^{\lambda}\lambda_i$.

In the following, we will use the notation $\mathcal{P}_+$ to denote
the set of all the partitions of positive integers. Let $0$ be the
partition of $0$, i.e. the empty partition. Define
$\mathcal{P}=\mathcal{P}_+\cup \{0\}$.

The power sum symmetric function of infinite variables
$\mathbf{x}=(x_1,..,x_N,..)$ is defined by
$p_{n}(\mathbf{x})=\sum_{i}x_i^n. $ Given a partition $\lambda$, we
define
$p_\lambda(\mathbf{x})=\prod_{j=1}^{l(\lambda)}p_{\lambda_j}(\mathbf{x}).
$ The Schur function $s_{\lambda}(\mathbf{x})$ is determined by the
Frobenius formula
\begin{align*}
s_\lambda(\mathbf{x})=\sum_{\mu}\frac{\chi_{\lambda}(\mu)}{\mathfrak{z}_\mu}p_\mu(\mathbf{x}),
\end{align*}
where $\chi_\lambda$ is the character of the irreducible
representation of
the symmetric group $S_{|\lambda|}$ corresponding to $\lambda$, we have $%
\chi_{\lambda}(\mu)=0$ if $|\mu|\neq |\lambda|$. The orthogonality
of character formula gives
\begin{align*}
\sum_\lambda\frac{\chi_\lambda(\mu)
\chi_\lambda(\nu)}{\mathfrak{z}_\mu}=\delta_{\mu \nu}.
\end{align*}
Let $n\in \mathbb{Z}$ and $\lambda,\mu,\nu$ denote the partitions.
We introduce the following notations
\begin{align*}
\{n\}_x=x^{\frac{n}{2}}-x^{-\frac{n}{2}}, \
\{\mu\}_{x}=\prod_{i=1}^{l(\mu)}\{\mu_i\}_x.
\end{align*}
In particular, let $\{n\}=\{n\}_q$ and $\{\mu\}=\{\mu\}_q$.

\subsection{LMOV invariants for framed knots}
Although the LMOV invariants are determined by the integrality
structure of topological open string partition function. Based on
the large $N$ duality of topological string and Chern-Simons theory
\cite{GV2,OV}, one can also introduce the LMOV invariants through
Chern-Simons theory of links/knots \cite{LM1,LM2,LMV,LP}.

Given a partition $\lambda$, we let
$\kappa_\lambda=\sum_{i=1}^{l(\lambda)}\lambda_i(\lambda_i-2i+1)$.
Let $\mathcal{K}_{\tau}$ be a knot with framing $\tau \in
\mathbb{Z}$. The framed colored HOMFLYPT invariant of
$\mathcal{K}_\tau$  is defined as follows
\begin{align*}
\mathcal{H}_\lambda(\mathcal{K}_\tau;q,a)=(-1)^{|\lambda|\tau}q^{\frac{\kappa_\lambda\tau}{2}}W_{\lambda}(\mathcal{K}_\tau;q,a),
\end{align*}
where $W_{\lambda}(\mathcal{K}_\tau;q,a)$ is the ordinary
(framing-independent) colored HOMFLYPT invariant of
$\mathcal{K}_\tau$, we refer to \cite{Zhu1} for the concrete
definition of $W_{\lambda}(\mathcal{K}_\tau;q,a)$.

Let $
\mathcal{Z}_{\mu}(\mathcal{K}_\tau;q,a)=\sum_{\lambda}\chi_{\lambda}(\mu)\mathcal{H}_{\lambda}(\mathcal{K}_\tau;q,a),
$ the Chern-Simons partition function of $\mathcal{K}_\tau$ is
defined by
\begin{align} \label{CSpartition}
Z_{CS}^{(S^3,\mathcal{K}_\tau)}(q,a,\mathbf{x})=\sum_{\lambda\in
\mathcal{P}}\mathcal{H}_\lambda(\mathcal{K}_\tau;q,a)s_{\lambda}(\mathbf{x})
=\sum_{\mu\in
\mathcal{P}}\frac{\mathcal{Z}_{\mu}(\mathcal{K}_\tau;q,a)}{\mathfrak{z}_\mu}p_{\mu}(\mathbf{x}).
\end{align}
 Then we define the functions
$f_{\lambda}(\mathcal{K}_\tau;q,a)$ by
\begin{align*}
Z_{CS}^{(S^3,\mathcal{K}_\tau)}(q,a,\mathbf{x})=\exp\left(\sum_{d=1}^\infty\frac{1}{d}\sum_{\lambda\in
\mathcal{P}_+}f_{\lambda}(\mathcal{K}_\tau;q^d,a^d)s_{\lambda}(\mathbf{x}^d)\right).
\end{align*}
Let $
\hat{f}_{\mu}(\mathcal{K}_\tau;q,a)=\sum_{\lambda}f_{\lambda}(\mathcal{K}_\tau;q,a)M_{\lambda\mu}(q)^{-1},
$ where
\begin{align*}
M_{\lambda\mu}(q)=\sum_{\mu}\frac{\chi_{\lambda}(C_{\nu})\chi_{\mu}(C_{\nu})}{\frak{z}_{\nu}}\prod_{j=1}^{l(\nu)}(q^{\nu_{j}/2}-q^{-\nu_{j}/2}).
\end{align*}

Denote $z=q^{\frac{1}{2}}-q^{-\frac{1}{2}}$, then the LMOV
conjecture for framed knot $\mathcal{K}_\tau$  stated that
\cite{MV}, for any $\mu\in \mathcal{P}_+$, there are integers
$N_{\mu,g,Q}(\mathcal{K}_\tau)$ such that
\begin{align*}
\hat{f}_{\mu}(\mathcal{K}_\tau;q,a)=\sum_{g\geq 0}\sum_{Q\in
\mathbb{Z}/2}N_{\mu,g,Q}(\mathcal{K}_\tau)z^{2g-2}a^Q\in
z^{-2}\mathbb{Z}[z^{2},a^{\pm \frac{1}{2}}].
\end{align*}
Therefore,
\begin{align*}
g_{\mu}(\mathcal{K}_\tau;q,a)&=\sum_{\nu}\chi_{\nu}(\mu)\hat{f}_\nu(\mathcal{K}_\tau;q,a)\\\nonumber
&=\sum_{g\geq 0}\sum_{Q\in
\mathbb{Z}/2}n_{\mu,g,Q}(\mathcal{K}_\tau)z^{2g-2}a^Q\in
z^{-2}\mathbb{Z}[z^{2},a^{\pm \frac{1}{2}}].
\end{align*}
where
$n_{\mu,g,Q}(\mathcal{K}_\tau)=\sum_{\nu}\chi_{\nu}(\mu)N_{\nu,g,Q}(\mathcal{K}_\tau)$.
These conjectural integers $n_{\mu,g,Q}(\mathcal{K}_\tau)$ ( and
$N_{\nu,g,Q}(\mathcal{K}_\tau)$) are referred to as the LMOV
invariants in this paper. We refer to \cite{GKS,KS1,KRSS1,KRSS2} for
another slightly different introduction of LMOV invariants which are
referred to as OV invariants in \cite{Zhu2}. The integrality of
certain OV invariants for a large family of knots/links have been
proved in \cite{KRSS1,KRSS2} recently by using the knots-quivers
correspondence.

In order to get an explicit expression for
$g_{\mu}(\mathcal{K}_\tau;q,a)$, we introduce
$F_{\mu}(\mathcal{K}_\tau;q,a)$ through the following expansion
formula
\begin{align*}
\log(Z_{CS}^{(S^3,\mathcal{K}_\tau)}(q,a,\mathbf{x}))=\sum_{\mu\in
\mathcal{P}_+}F_{\mu}(\mathcal{K}_\tau;q,a)p_{\mu}(\mathbf{x}).
\end{align*}
Then, by formula (\ref{CSpartition}), we have
\begin{align*}
F_{\mu}(\mathcal{K}_\tau;q,a)=\sum_{n\geq
1}\sum_{\cup_{i=1}^{n}\nu^i=\mu}\frac{(-1)^{n-1}}{n}\prod_{i=1}^{n}\frac{\mathcal{Z}_{\nu^i}(\mathcal{K}_\tau;q,a)}{\mathfrak{z}_{\nu^i}}.
\end{align*}

\begin{remark}
For two partitions $\nu^1$ and $\nu^2$, the notation $\nu^1\cup
\nu^2$ denotes the new partition obtained by combining all the parts
in $\nu^1, \nu^2$. For example, if $\mu=(2,2,1)$, then the list of
all the pairs $(\nu^1,\nu^2)$ such that $\nu^1 \cup \nu^2 =(2,2,1)$
is
\begin{align*}
(\nu^1=(2), \nu^2=(2,1)),\ (\nu^1=(2,1), \nu^2=(2)),\\\nonumber
 \
(\nu^1=(1), \nu^2=(2,2)), \ (\nu^1=(2,2), \nu^2=(1)). \
\end{align*}
\end{remark}

Finally, by using the M\"obius inversion formula, we obtain
\begin{align} \label{formulagmu}
g_{\mu}(\mathcal{K}_\tau;q,a)=\mathfrak{z}_{\mu}\frac{1}{\{\mu\}}\sum_{d|\mu}\frac{\mu(d)}{d}F_{\mu/d}(\mathcal{K}_\tau;q^d,a^d).
\end{align}
\begin{remark}
In the above discussion, for the sake of brevity, we only consider
the case of a framed knot, actually, the LMOV invariants can be
defined for any framed link.
\end{remark}

\subsection{LMOV invariants for framed unknot $U_\tau$}
In the following, we only consider the case of a framed unknot
$U_\tau$. In this situation, the large $N$ duality of topological
string and Chern-Simons theory \cite{MV} has been proved in
\cite{LLZ,Zhou}. Therefore, we can also compute the LMOV invariants
for framed unknot $U_\tau$ through the open topological string
theory.

We denote by $n_{m,l}(\tau)$ the LMOV invariants
$n_{(m),0,l-\frac{m}{2}}(U_\tau)$ of the framed unknot $U_\tau$,
where $g=0$ and $m\geq 1$, $l\geq 0$. According to the computations
shown in \cite{MV} ( or cf. pages 15-16 in \cite{LZ}), the explicit
closed formula for $n_{m,l}(\tau)$ is given by formula
(\ref{integralformula}). The integrality of $n_{m,l}(\tau)$ is given
by Theorem \ref{Thm1}.

Then, denoted by $n_{(m_1,m_2)}(\tau)$ the LMOV invariants
$n_{(m_1,m_2),0,\frac{m_1+m_2}{2}}(U_\tau)$ of the framed unknot
$U_\tau$ with  $\mu=(m_1,m_2)$, $g=0$ and $Q=\frac{m_1+m_2}{2}$,
where $m_1\geq m_2\geq 1$. By using the computations in open
topological string theory (cf. pages 18-19 in \cite{LZ}), we obtain
the explicit closed formula for $n_{(m_1,m_2)}(\tau)$ which is given
by formula (\ref{integralformula2}).   The integrality of
$n_{(m_1,m_2)}(\tau)$ is given by Theorem \ref{Thm2}.

Let $n_{m,g,Q}(\tau)$ be the LMOV invariants $n_{(m),g,Q}(U_\tau)$
of the framed unknot $U_\tau$. We consider the following generating
function
\begin{align*}
g_m(q,a)=\sum_{g\geq 0}\sum_{Q\in
\mathbb{Z}/2}n_{m,g,Q}(\tau)z^{2g-2}a^Q
\end{align*}
which can be computed by using formula (\ref{formulagmu}).
Considering the following function
\begin{align*}
\phi_{\mu,\nu}(x)=\sum_{\lambda}\chi_{\lambda}(\mu)\chi_{\lambda}(\nu)x^{\kappa_\lambda}.
\end{align*}
By Lemma 5.1  in \cite{CLPZ},  for $d\in \mathbb{Z}_+$, we have
\begin{align} \label{formulaphi}
\phi_{(d),\nu}(x)=\frac{\{d\nu\}_{x^2}}{\{d\}_{x^2}}.
\end{align}
By the expression of colored HOMFLYPT invariant for unknot (cf.
formula (4.6) in \cite{LP})
\begin{align*}
W_{\lambda}(U;q,a)=\sum_{\nu}\frac{\chi_{\lambda}(\nu)}{\mathfrak{z}_\nu}\frac{\{\nu\}_a}{\{\nu\}},
\end{align*}
we obtain
\begin{align*}
\mathcal{Z}_{\mu}(U_\tau;q,a)&=\sum_{\lambda}\chi_{\lambda}(\mu)\mathcal{H}_{\lambda}(U_\tau;q,a)\\\nonumber
&=(-1)^{|\mu|\tau}\sum_{\lambda}\chi_{\lambda}(\mu)q^{\frac{\kappa_\lambda
\tau}{2}}\sum_{\nu}\frac{\chi_{\lambda}(\nu)}{\mathfrak{z}_\nu}\frac{\{\nu\}_a}{\{\nu\}}\\\nonumber
&=(-1)^{|\mu|\tau}\sum_{\nu}\frac{1}{\mathfrak{z}_\nu}\phi_{\mu,\nu}(q^\frac{\tau}{2})\frac{\{\nu\}_a}{\{\nu\}}.
\end{align*}
In particular, for $\mu=(m)$ with $m\geq 1$, formula
(\ref{formulaphi}) implies that
\begin{align*}
\mathcal{Z}_{(m)}(U_\tau;q,a)=(-1)^{m\tau}\sum_{|\nu|=m}\frac{1}{\mathfrak{z}_\nu}\frac{\{m\nu\tau\}}{\{m\tau\}}\frac{\{\nu\}_a}{\{\nu\}}.
\end{align*}
For brevity, if we let
\begin{align*}
\mathcal{Z}_m(q,a)=\frac{1}{\{m\}}\mathcal{Z}_{(m)}(U_\tau;q,a)=(-1)^{m\tau}\sum_{|\nu|=m}\frac{1}{\mathfrak{z}_\nu}
\frac{\{m\nu\tau\}}{\{m\}\{m\tau\}}\frac{\{\nu\}_a}{\{\nu\}},
\end{align*}
then formula (\ref{formulagmu}) gives
\begin{align} \label{gm}
g_m(q,a)=\sum_{d|m}\mu(d)\mathcal{Z}_{m/d}(q^d,a^d).
\end{align}
Theorem \ref{Thm3} shows that
 $ g_m(q,a)\in z^{-2}\mathbb{Z}[z^2,a^{\pm \frac{1}{2}}]$ for any $m\geq
 1$.

\section{Proof of the Theorem \ref{Thm1}} \label{section3}
For nonnegative integer $n$ and prime number $p$, we introduce the
following function
\begin{equation} \label{functionfp}
    f_p(n)=\prod_{i=1,p\nmid i}^n i = \frac{n!}{p^{[n/p]}[n/p]!}.
\end{equation}
Given a positive integer $k$, throughout this paper, we use the
notation $p^k \mid\mid n$ to denote that $p^k$ divides $n$, but
$p^{k+1}$ does not.

Before giving the proof of Theorem \ref{Thm1}, we first establish
several useful lemmas.

\begin{lemma}
    For odd prime numbers $p$ and $\alpha\geq 1$ or for $p=2$, $\alpha\geq 2$,
    we have $p^{2\alpha}\mid f_p(p^{\alpha} n)-f_p(p^{\alpha})^n$. For $p=2, \alpha=1$, $f_2(2n)\equiv (-1)^{[n/2]}\pmod{4}$
    \label{lemma1}.
\end{lemma}

\begin{proof}
 With $\alpha\geq 2$ or $p>2$, $p^{\alpha-1}(p-1)$ is even,
 \begin{align*}
     &f_p(p^\alpha n)-f_p(p^\alpha(n-1))f_p(p^\alpha) \\
     &= f_p(p^\alpha(n-1))\left(\prod_{i=1,p\nmid i}^{p^\alpha} (p^{\alpha}(n-1)+i)-f_p(p^\alpha)\right) \\
     & \equiv p^{\alpha}(n-1) f_p(p^\alpha(n-1))f_p(p^\alpha)\left(\sum_{i=1,p\nmid i}^{p^\alpha}\frac{1}{i}\right) \pmod{p^{2\alpha}} \\
     & \equiv  p^{\alpha}(n-1) f_p(p^\alpha(n-1))f_p(p^\alpha)\left(\sum_{i=1,p\nmid i}^{[p^\alpha/2]} (\frac{1}{i}+\frac{1}{p^{\alpha}-i})\right) \pmod{p^{2\alpha}} \\
     & \equiv  p^{\alpha}(n-1) f_p(p^\alpha(n-1))f_p(p^\alpha)\left(\sum_{i=1,p\nmid i}^{[p^\alpha/2]} \frac{p^{\alpha}}{i(p^\alpha-i)} \right)\equiv 0, \pmod{p^{2\alpha}}
\end{align*}

Thus the first part of the Lemma \ref{lemma1} is proved by
induction. For $p=2, \alpha=1$, the formula is straightforward.
\end{proof}

\begin{lemma}
    For odd prime number $p$ and $m=p^\alpha a, l=p^\beta b$, $p\nmid a, p\nmid b$, $\alpha\geq 1, \beta \geq 0$, we have
    \begin{equation*}
        p^{2\alpha} \mid \binom{m}{l}\binom{m\tau+l-1}{m-1}-\binom{\frac{m}{p}}{\frac{l}{p}}\binom{\frac{m\tau+l}{p}-1}{\frac{m}{p}-1}
    \end{equation*}
    where for $\beta=0$, the second term is defined to be zero. \label{luo:lemma2-0}
\end{lemma}
\begin{proof}
    \begin{align}
        & \binom{m}{l}\binom{m\tau+l-1}{m-1}-\binom{\frac{m}{p}}{\frac{l}{p}}\binom{\frac{m\tau+l}{p}-1}{\frac{m}{p}-1} \nonumber \\
        & = \binom{\frac{m}{p}}{\frac{l}{p}}\binom{\frac{m\tau+l}{p}-1}{\frac{m}{p}-1} \left( \frac{f_p(m)}{f_p(l)f_p(m-l)}\cdot \frac{f_p(m\tau+l)}{f_p(m)f_p(m(\tau-1)+l)}-1\right) \label{luo:eq5-0}
    \end{align}

    Write $\binom{\frac{m}{p}}{\frac{l}{p}}=\frac{m}{l}\binom{\frac{m}{p}-1}{\frac{l}{p}-1}$ and $\binom{\frac{m\tau+l}{p}-1}{\frac{m}{p}-1}=\frac{m}{m\tau +l} \binom{\frac{m\tau+l}{p}}{\frac{m}{p}}$, both are divisible by $p^{\max(\alpha-\beta,0)}$.
    Each element of $\{m, l, m-l, m\tau+l, m(\tau-1)+l\}$ is divisible by $p^{\min(\alpha,\beta)}$, so by Lemma~\ref{lemma1},
\begin{equation}
    \frac{f_p(m)}{f_p(l)f_p(m-l)}\cdot \frac{f_p(m\tau+l)}{f_p(m)f_p(m(\tau-1)+l)}-1 \label{luo:term1-0}
\end{equation}
is divisible by $p^{2\min(\alpha,\beta)}$ (including the case
$\beta=0$) in $p$-adic number field. Thus (\ref{luo:eq5-0}) is
divisible by
$p^{2\max(\alpha-\beta,0)+2\min(\alpha,\beta)}=p^{2\alpha}$.
\end{proof}

\begin{lemma}
    For $m=2^\alpha a, l=2^\beta b, \alpha\geq 1, \beta \geq 0$,
    \[ 2^{2\alpha}\mid (-1)^{m\tau+m+l}\binom{m}{l}\binom{m\tau+l-1}{m-1}-(-1)^{\frac{m\tau+m+l}{2}}\binom{\frac{m}{2}}{\frac{l}{2}}\binom{\frac{m\tau+l}{2}-1}{\frac{m}{2}-1},\]
    where the second term is set to zero for $\beta=0$. \label{luo:lemma3-0}
\end{lemma}

\begin{proof}
    For the case $\alpha\geq 2, \beta\geq 2$, both $m\tau+m+l$ and $(m\tau+m+l)/2$ are even, the Lemma is proved as in Lemma~\ref{luo:lemma2-0}. For the case $\beta=0$, both $\binom{m}{l}$ and $\binom{m\tau+l-1}{m-1}$ are divisible by $2^\alpha$, and the Lemma is also proved.
    For remaining cases $\alpha>\beta =1$ or $\beta\geq \alpha=1$, we compute similarly as (\ref{luo:eq5-0}),
    \begin{align}
        & (-1)^{m\tau+m+l}\binom{m}{l}\binom{m\tau+l-1}{m-1}-(-1)^{\frac{m\tau+m+l}{2}}\binom{\frac{m}{2}}{\frac{l}{2}}\binom{\frac{m\tau+l}{2}-1}{\frac{m}{2}-1}  \nonumber \\
        & = \binom{\frac{m}{2}}{\frac{l}{2}}\binom{\frac{m\tau+l}{2}-1}{\frac{m}{2}-1} \left( \frac{f_2(m)}{f_2(l)f_2(m-l)}\cdot \frac{f_2(m\tau+l)}{f_2(m(\tau-1)+l)f_2(m)}-(-1)^{\frac{m\tau+m+l}{2}}\right) \label{luo:eq6-0}
    \end{align}
    Both $\binom{\frac{m}{2}}{\frac{l}{2}}$ and $\binom{\frac{m\tau+l}{2}-1}{\frac{m}{2}-1}$ are divisible by $2^{\alpha-1}$, it suffices to prove that the third factor is divisible by $4$, which is, by Lemma~\ref{lemma1},
    \begin{equation}
        (-1)^{[\frac{l}{4}]+[\frac{m-l}{4}]+[\frac{m\tau+l}{4}]+[\frac{m(\tau-1)+l}{4}]}-(-1)^{\frac{m\tau+m+l}{2}}. \pmod{4} \nonumber
    \end{equation}
    It is divisible by $4$ if
    \begin{equation}
        [\frac{l}{4}]+[\frac{m-l}{4}]+[\frac{m\tau+l}{4}]+[\frac{m(\tau-1)+l}{4}]+\frac{m\tau+m+l}{2}\label{luo:eq7-0}
    \end{equation}
    is even.
Parity of (\ref{luo:eq7-0})  depends only on $\tau\pmod{2}$. For
$\tau=1$, (\ref{luo:eq7-0}) reduces to
$[l/4]+[(m-l)/4]+[(m+l)/4]+[l/4]+l/2$. For $\tau=0$, it reduces to
$[l/4]+[(m-l)/4]+[l/4]+[(l-m)/4]+(m+l)/2$. Both are obviously even.

\end{proof}
Now, we can finish the proof of Theorem \ref{Thm1}.

\begin{proof}
    For a prime number $p\mid m$, write $m=p^\alpha a, p\nmid a$.
    \begin{align}
        n_{m,l}(\tau)&=\sum_{d\mid m, d\mid l} \frac{\mu(d)}{d^2} c_{\frac{m}{d},\frac{l}{d}}(\tau)  \nonumber \\
                     &= \frac{1}{m^2} \sum_{d\mid m, d\mid l}\mu(d)(-1)^{\frac{m\tau+m+l}{d}} \binom{\frac{m}{d}}{\frac{l}{d}}\binom{\frac{m\tau+l}{d}-1}{\frac{m}{d}-1} \nonumber \\
        &=\frac{1}{m^2} \sum_{d\mid m, d\mid l, p\nmid d} \mu(d) \left((-1)^{\frac{m\tau+m+l}{d}}\binom{\frac{m}{d}}{\frac{l}{d}}\binom{\frac{m\tau+l}{d}-1}{\frac{m}{d}-1}-(-1)^{\frac{m\tau+m+l}{dp}} \binom{\frac{m}{dp}}{\frac{l}{dp}}\binom{\frac{m\tau+l}{dp}-1}{\frac{m}{dp}-1}\right) \label{luo:eq9-0}
    \end{align}
    where for $pd\nmid l$, the second term of (\ref{luo:eq9-0}) is understood to be zero. For odd prime number $p$, $\frac{m\tau+m+l}{d}$ and $\frac{m\tau+m+l}{dp}$ have the same parity. Since $p^\alpha\mid\mid\frac{m}{d}$, $p^{2\alpha}$ divides the summand in (\ref{luo:eq9-0})  by Lemma~\ref{luo:lemma2-0}. For $p=2$, it is divisible by $2^{2\alpha}$ by  Lemma~\ref{luo:lemma3-0}.

    $p^{2\alpha}$ divides the sum in (\ref{luo:eq9-0}) for every $p^\alpha \mid\mid m$, thus $n_{m,l}$ is an integer.
\end{proof}

\section{Proof of the Theorem \ref{Thm2}} \label{section4}

We introduce the following lemma first.
\begin{lemma}
    If $p^\beta\mid\mid (a,b), p^\alpha\mid a+b$, then $p^{\alpha-\beta}$ divides
    $$ \binom{a\tau+a-1}{a}\binom{b\tau+b}{b}.$$
    \label{luo:lemma3-1}
\end{lemma}
\begin{proof}
    Power of prime $p$ in $n!$ is
    $ \sum_{k=1}^\infty [\frac{n}{p^k}]. $
    Apply this to the binomial coefficients to find that the power of $p$ in $\binom{a\tau+a-1}{a}\binom{b\tau+b}{b}$ is
    \begin{align*}
        & \sum_{i=1}^{\infty} \left([\frac{a\tau+a-1}{p^i}]+[\frac{b\tau+b}{p^i}]\right)-\left([\frac{a\tau-1}{p^i}]+[\frac{b\tau}{p^i}]\right)-\left([\frac{a}{p^i}]+[\frac{b}{p^i}]\right) \\
        &\geq \sum_{i=1}^{\alpha} \left((\frac{(a+b)(\tau+1)}{p^i}-1)-(\frac{(a+b)\tau}{p^i}-1)\right) - \sum_{i=1}^{\beta}(\frac{a+b}{p^i})-\sum_{i=\beta+1}^{\alpha}(\frac{a+b}{p^i}-1) \\
        &=\alpha-\beta
    \end{align*}
    where we use the fact that for $k\mid m+n+1, k>1$, $[m/k]+[n/k]=(m+n+1)/k-1$ and for $k\mid m+n, k\nmid m$, $[m/k]+[n/k]=(m+n)/k-1$.
\end{proof}

Recall the definition of the function $f_{p}(n)$ given by formula
(\ref{functionfp}). It is obvious that
\begin{equation}
    f_p(p^\alpha k) \equiv f_p(p^{\alpha})^k \equiv (-1)^k \pmod{p^\alpha} \label{luo:eq3-1}
\end{equation}

Now, we can finish the proof of  Theorem \ref{Thm2}:
\begin{proof}
 By definition,
\begin{align}
    n_{(m_1,m_2)}(\tau)& =  \frac{1}{m_1+m_2}\sum_{d\mid m_1,d\mid m_2} \mu(d) (-1)^{(m_1+m_2)(\tau+1)/d} \nonumber \\
                     & \cdot \binom{(m_1\tau+m_1)/d-1}{m_1/d}\binom{(m_2\tau+m_2)/d}{m_2/d} \label{luo:eq1-1}
\end{align}

    Let $p$ be any prime divisor of $m_1+m_2$, $p^\alpha\mid \mid m_1+m_2$. We will prove $p^\alpha$ divides the summation in (\ref{luo:eq1-1}), thus $m_1+m_2$ also divides and $n_{m_1,m_2}$ are integers.

   If $p\nmid m_1$, each summand in (\ref{luo:eq1-1}) corresponds to $p\nmid d$, so $p^\alpha\mid (m_1+m_2)/d$ and $p\nmid m_1/d$. By Lemma~\ref{luo:lemma3-1} applied to $a=m_1/d, b=m_2/d$, $p^{\alpha}$ divides each summand and thus the summation.

   If $p^\beta \mid\mid m_1, \beta\geq 1$, consider two summands in (\ref{luo:eq1-1}) corresponding to $d$ and $pd$ such that $pd\mid (m_1,m_2), \mu(pd)\neq 0$. When $p$ is an odd prime or $\alpha\geq 2$, the sign $(-1)^{(m_1+m_2)(\tau+1)/d}$ and $(-1)^{(m_1+m_2)(\tau+1)/(pd)}$ are equal. When $p=2, \alpha=1$, modulo 2 the sign is irrelevant. Write $a=m_1/d, b=m_2/d$, then $p^\alpha\mid a+b, p^\beta\mid\mid a$.
   \begin{align}
       & \binom{a\tau+a-1}{a}\binom{b\tau+b}{b}-\binom{(a\tau+a)/p-1}{a/p}\binom{(b\tau+b)/p}{b/p} \nonumber \\
       &= \binom{(a\tau+a)/p-1}{a/p}\binom{(b\tau+b)/p}{b/p}\left( \frac{f_p(a\tau+a)f_p(b\tau+b)}{f_p(a\tau)f_p(a)f_p(b\tau)f_p(b)} - 1\right)  \nonumber \\
       &= \binom{(a\tau+a)/p-1}{a/p}\binom{(b\tau+b)/p}{b/p} \frac{f_p(a\tau+a)f_p(b\tau+b)-f_p(a\tau)f_p(a)f_p(b\tau)f_p(b)}{f_p(a\tau)f_p(a)f_p(b\tau)f_p(b)}  \label{luo:eq2-1}
   \end{align}

   The term $\binom{(a\tau+a)/p-1}{a/p}\binom{(b\tau+b)/p}{b/p}$ is divisible by $p^{\alpha-\beta}$ by Lemma~\ref{luo:lemma3-1}. The numerator of the fraction term in (\ref{luo:eq2-1}) is divisible by $p^\beta$ by (\ref{luo:eq3-1}), and the denominator is not divisible by $p$. We proved that $p^\alpha$ divides (\ref{luo:eq2-1}), take summation over $d$, we get that $p^\alpha$ divides the summation in (\ref{luo:eq1-1}).
    This is true for any $p\mid m_1+m_2$, thus $n_{(m_1,m_2)}(\tau)$ is an integer.
\end{proof}

\section{Proof of the Theorem \ref{Thm3}} \label{section5}
We establish several lemmas first.
\begin{lemma}
    Suppose $k$ is a positive integer, then the number
    \begin{equation*}
        c_m(k,y)=\sum_{|\lambda|=m} \frac{1}{\mathfrak{z}_\lambda} k^{l(\lambda)}\{\lambda\}_{y^2}
    \end{equation*}
    is equal to the coefficient of $t^m$ in $(\frac{1-t/y}{1-ty})^k$.
\end{lemma}

\begin{proof}
    Suppose the number of $i$'s in the partition $\lambda$ is $a_i, i=1,\cdots$. Then
    \begin{align*}
        c_m(k,y)&=\sum_{\sum ia_i=m} \prod_i \frac{1}{a_i!i^{a_i}} k^{a_i}(y^i-y^{-i})^{a_i} \\
                &=\left[ \prod_{i=1}^{\infty} \left(\sum_{j=0}^\infty t^{ij}\frac{1}{j!i^j} k^j(y^i-y^{-i})^j\right) \right]_{t^m} \\
                &=\left[ \prod_{i=1}^{\infty} \exp(t^ik(y^i-y^{-i})/i) \right]_{t^m} \\
                &=\left[ \exp(k\ln(1-ty)^{-1}+k\ln(1-t/y))\right]_{t^m}\\
                &=\left[ (\frac{1-t/y}{1-ty})^k \right]_{t^m}
    \end{align*}
\end{proof}

\begin{lemma} \label{mmtaulemma}
Let $R=\mathbb{Q}[q^{\pm 1/2},a^{\pm 1/2}]$. Then
\begin{align}  \label{mmtaugm}
 \{m\}\{m\tau\}g_m(q,a)=\sum_{d\mid m}\sum_{|\mu|=m/d}\frac{\mu(d)(-1)^{m\tau/d}}{\mathfrak{z}_\mu}\frac{\{m\mu\tau\}}{\{d\mu\}}\{d\mu\}_a
\end{align}
is divisible by $\{m\tau\}\{m\}/\{1\}^2$ in $R$.
\end{lemma}
\begin{proof}
By the definition (\ref{gm}) of $g_m(q,a)$, we have the formula
(\ref{mmtaugm}). It is clear that
$$\{m\}\{m\tau\}g_m(q,a) \in R.$$
Denote $\Phi_n(q)=\prod_{d\mid n} (q^d-1)^{\mu(n/d)}$ to be the
$n$-th cyclotomic polynomial, which is irreducible over $R$. Then
$q^n-1=\prod_{d\mid n} \Phi_d(q)$, and
\begin{align} \label{cyc}
   \{m\}\{m\tau\}&=q^{-\frac{m+m\tau}{2}} \prod_{{m_1}\mid m} \Phi_{m_1}(q)\prod_{{m_1}\mid m\tau} \Phi_{m_1}(q)\\\nonumber
   &=q^{-\frac{m+m\tau}{2}}\prod_{{m_1}\mid m} \Phi_{m_1}(q)^2\prod_{{m_1}\mid m\tau,
{m_1}\nmid m} \Phi_{m_1}(q)
\end{align}

(i) For $m_1\mid m\tau, m_1\nmid m$, and any $|\mu|=m/d$, at least
one of $d\mu_i$'s are not divisible by $m_1$, thus
$\{m\mu_i\tau\}/\{d\mu_i\}$ is divisible by $\Phi_{m_1}(q)$. So
$\Phi_{m_1}(q)$ divides $\{m\}\{m\tau\}g_m(q,a)$.

(ii) For $m_1\mid m$ and any $|\mu|=m/d$, if not all $d\mu_i$ are
divisible by $m_1$, then at least two of them are not divisible.
Then two of corresponding $\{m\mu_i\tau\}/\{d\mu_i\}$ are divisible
by $\Phi_{m_1}(q)$.

We consider modulo $\{m_1\}^2$ in the ring $R$. It is easy to see,
for $a,b\geq 1 $,
\begin{align*}
\frac{\{abm_1\}}{\{bm_1\}} \equiv
a\left(\frac{q^{m_1/2}+q^{-m_1/2}}{2}\right)^{(a-1)b}
\pmod{\{m_1\}^2}
\end{align*}
We write $x=(q^{m_1/2}+q^{-m_1/2})/2$, then $x^2\equiv 1
\pmod{\{m_1\}^2}$.

Then modulo $\Phi_{m_1}(q)^2$, we have
\begin{align}
    &\{m\}\{m\tau\}g_m(q,a)\nonumber\\
     &\equiv \sum_{d\mid m}\sum_{|\mu|=m/d,m_1\mid d\mu} \frac{\mu(d)(-1)^{m\tau/d}}{\mathfrak{z}_\mu} \frac{\{m\mu\tau\}}{\{d\mu\}}\,\{d\mu\}_a \nonumber\\
                          &\equiv \sum_{d\mid m}\sum_{|\mu|=m/d,m_1\mid d\mu} \frac{\mu(d)(-1)^{m\tau/d}}{\mathfrak{z}_\mu} \left(\frac{m\tau}{d}\right)^{l(\mu)}
                          x^{(m|\mu|\tau-d|\mu|)/m_1} \{d\mu\}_a  \nonumber\\
                          &\equiv \sum_{d\mid m}\sum_{|\lambda|=m/\mathrm{lcm}(d,m_1)} \mu(d)(-1)^{m\tau/d}x^{\frac{m}{m_1}(\frac{m\tau}{d}-1)}
                          \cdot\frac{1}{\mathfrak{z}_\lambda}\left(\frac{m\tau}{\mathrm{lcm}(d,m_1)}\right)^{l(\lambda)} \{\lambda\}_{a^{\mathrm{lcm}(d,m_1)}} \nonumber\\
                          &\equiv \sum_{d\mid m} \mu(d) (-1)^{m\tau/d} x^{\frac{m}{m_1}(\frac{m\tau}{d}-1)} \left[ \left(\frac{1-t^{\mathrm{lcm}(d,m_1)}a^{-\mathrm{lcm}(d,m_1)/2}}{1-t^{\mathrm{lcm}(d,m_1)}a^{\mathrm{lcm}(d,m_1)/2}}\right)^{m\tau/\mathrm{lcm}(d,m_1)} \right]_{t^m}
                          \label{luo:eq44}
\end{align}
\begin{itemize}
    \item For the cases $m_1$ with an odd prime factor $p$, or $p=2$ divides $m_1$ and $4\mid m$, or $p=2$ divides $m_1$ and $2\mid \tau$: Consider those $d$ with $\mu(d)\neq 0$ and $p\nmid d$, we have $\mathrm{lcm}(d,m_1)=\mathrm{lcm}(pd,m_1)$ and parity of $m\tau/d$ equals parity of $m\tau/(pd)$, but $\mu(d)=-\mu(pd)$. Thus two terms in (\ref{luo:eq44}) corresponding to $d$ and $pd$ cancelled.
    \item For the remaining case $2\mid\mid m, m_1=2, 2\nmid\tau$: $\Phi_{m_1}(q)^2=(q^{1/2}+q^{-1/2})^2=2x+2$. Coefficients of $x$ in (\ref{luo:eq44}) equals sum of terms corresponds to odd $d\mid m, \mu(d)\neq 0$, while constant term coefficient equals to sum of terms corresponds to $2d\mid m, \mu(2d)\neq 0$. The coefficients of term for $d$ and $2d$ match, so (\ref{luo:eq44}) is divisible by $x+1$.
\end{itemize}

In summary, we have proved that for $m_1\mid m\tau, m_1\nmid m$,
$\Phi_{m_1}(q)$ divides $\{m\}\{m\tau\}g_m(q,a)$; for $m_1\mid m,
m_1 \neq 1$, $\Phi_{m_1}(q)^2$ divides $\{m\}\{m\tau\}g_m(q,a)$. By
(\ref{cyc}), the lemma is proved.
\end{proof}

\begin{lemma} \label{lemmaprime}
For any integer $m\geq 1$, we have
\begin{align*}
g_m(q,a)\in z^{-2}\mathbb{Q}[z^2,a^{\pm \frac{1}{2}}].
\end{align*}
\end{lemma}
\begin{proof}
By Lemma \ref{mmtaulemma}, we have
\begin{align*}
f(q,a):=z^2g_m(q,a)=\frac{\{1\}^2}{\{m\}\{m\tau\}}\sum_{d\mid
 m}\sum_{|\mu|=m/d}\frac{\mu(d)(-1)^{m\tau/d}}{\mathfrak{z}_\mu}\frac{\{m\mu\tau\}}{\{d\mu\}}\{d\mu\}_a
 \in \mathbb{Q}[q^{\pm \frac{1}{2}},a^{\pm 1}].
\end{align*}
As a function of $q$, it is clear $f(q,a)$ admits $
f(q,a)=f(q^{-1},a). $ Furthermore, for any $d|m$ and $|\mu|=m/d$, we
have
\begin{align*}
   m|\mu|\tau-d|\mu|-m\tau-m\equiv m^2\tau/d-m\tau=m\tau(m/d-1) \equiv 0
   \pmod{2},
\end{align*}
which implies $ f(q,a)=f(-q,a). $ Therefore,
$f(q,a)=z^2g_{m}(q,a)\in \mathbb{Q}[z^2,a^{\pm \frac{1}{2}}]$. The
lemma is proved.
\end{proof}

\begin{lemma} \label{integralitylemma}
For any $\tau\in \mathbb{Z}$, we have
\begin{align} \label{integralityformula}
\{m\}\{m\tau\}\mathcal{Z}_m(q,a)\in \mathbb{Z}[q^{\pm
\frac{1}{2}},a^{\pm \frac{1}{2}}].
\end{align}
\end{lemma}

\begin{proof}
Since
\begin{align*}
(-1)^{m\tau}\{m\}\{m\tau\}\mathcal{Z}_m(q,a)&=\sum_{|\mu|=m}\frac{\{m\tau
\mu\}}{\mathfrak{z}_\mu\{\mu\}}\{\mu\}_a\\\nonumber
&=\sum_{\sum_{j\geq 1}jk_j=m}\frac{\prod_{j\geq 1}(\{m\tau
j\}\{j\}_a)^{k_j}}{\prod_{j\geq 1}j^{k_j}k_j!},
\end{align*}
we construct a generating function
\begin{align} \label{formulafm}
f(x)&=\sum_{n\geq 0}x^n\sum_{\sum_{j\geq 1}jk_j=n}\frac{\prod_{j\geq
1}(\{m\tau j\}\{j\}_a)^{k_j}}{\prod_{j\geq 1}j^{k_j}k_j!}\\\nonumber
&=\sum_{n\geq 0}\sum_{\sum_{j\geq 1}jk_j=n}\frac{\prod_{j\geq
1}(\{m\tau j\}\{j\}_ax^j)^{k_j}}{\prod_{j\geq
1}j^{k_j}k_j!}\\\nonumber &=\exp\left(\sum_{j\geq 1}\frac{\{m\tau
j\}\{j\}_ax^j}{j\{j\}}\right),
\end{align}
Then $
(-1)^{m\tau}\{m\}\{m\tau\}\mathcal{Z}_m(q,a)=\left[f(x)\right]_{x^m}.
$

For $\tau=0$, it is the trivial case.

For $\tau \geq 1$, we use the expansion $\frac{\{m\tau
j\}}{\{j\}}=\sum_{k=0}^{m\tau-1}q^{\frac{j(m\tau-2k-1)}{2}}$, then
\begin{align*}
f(x)&=\exp\left(\sum_{k\geq 0}^{m\tau-1}\sum_{j\geq
1}\left(\frac{(q^{\frac{m\tau-1-2k}{2}}a^{\frac{1}{2}}x)^j}{j}-\frac{(q^{\frac{m\tau-1-2k}{2}}a^{-\frac{1}{2}}x)^j}{j}\right)\right)\\\nonumber
&=\exp\left(\sum_{k\geq
0}^{m\tau-1}\log\frac{1+q^{\frac{m\tau-1-2k}{2}}a^{-\frac{1}{2}}x}{1+q^{\frac{m\tau-1-2k}{2}}a^{\frac{1}{2}}x}\right)\\\nonumber
&=\prod_{k=0}^{m\tau-1}\frac{1+q^{\frac{m\tau-1-2k}{2}}a^{-\frac{1}{2}}x}{1+q^{\frac{m\tau-1-2k}{2}}a^{\frac{1}{2}}x}.
\end{align*}
We introduce the $q$-binomial coefficients defined by
\begin{align*}
\binom{m}{r}_q=\frac{(1-q^m)(1-q^{m-1})\cdots
(1-q^{m-r+1})}{(1-q)(1-q^2)\cdots (1-q^r)}
\end{align*}
for $r\leq m$, and in particular $\binom{m}{0}_q=1$. The
$q$-binomial coefficients $\binom{m}{r}_q\in \mathbb{Z}[q]$ (see
Chapter 2 of \cite{KS0} for $q$-binomial coefficients). There are
analogs of the binomial formula, and of Newton's generalized version
of it for negative integer exponents,
\begin{align*}
\prod_{k=0}^{n-1}(1+q^kt)&=\sum_{k=0}^nq^{\frac{k(k-1)}{2}}\binom{n}{k}_qt^k
\\\nonumber
\prod_{k=0}^{n-1}\frac{1}{(1-q^kt)}&=\sum_{k=0}^\infty
\binom{n+k-1}{k}_qt^k.
\end{align*}

Therefore, the coefficient $\left[f(x)\right]_{x^m}$ of $x^m$ in
$f(x)$ is given by
\begin{align*}
\sum_{j+k=m}(-1)^{k}q^{\frac{j(j-1)-(m\tau-1)m}{2}}a^{\frac{k-j}{2}}\binom{m\tau}{j}_q\binom{m\tau+k-1}{k}_q,
\end{align*}
which lies in the ring $\mathbb{Z}[q^{\pm \frac{1}{2}},a^{\pm
\frac{1}{2}}]$ by the integrality of Gaussian binomial.

For the case $\tau \leq -1$, we write $\{m\tau j\}=-\{-m\tau j\}$ in
the formula (\ref{formulafm}), then the similar computations give
the formula (\ref{integralityformula}).
\end{proof}

Now, we can finish the proof of Theorem \ref{Thm3} as follow:
\begin{proof}
Lemma \ref{lemmaprime}  implies that there exist
 rational numbers $n_{m,g,Q}(\tau)$, such that
\begin{align*}
z^2g_{m}(q,a)=\sum_{g \geq 0}\sum_{Q}n_{m,g,Q}(\tau)z^{2g}a^{Q}\in
\mathbb{Q}[z^2,a^{\pm \frac{1}{2}}].
\end{align*}
So we only need to show $n_{m,g,Q}(\tau)$ are integers. By lemma
\ref{integralitylemma} and the formula (\ref{gm}) for $g_{m}(q,a)$,
we have
\begin{align*}
\{m\}\{m\tau\}z^2g_{m}(q,a)\in \mathbb{Z}[q^{\pm \frac{1}{2}},a^{\pm
\frac{1}{2}}],
\end{align*}
which is equivalent to
\begin{align*}
(q^{\frac{m}{2}}-q^{-\frac{m}{2}})(q^{\frac{m\tau}{2}}-q^{-\frac{m\tau}{2}})\sum_{g
\geq 0}\sum_{Q}n_{m,g,Q}(\tau)(q^{1/2}-q^{-1/2})^{2g}a^{Q}\in
\mathbb{Z}[q^{\pm \frac{1}{2}},a^{\pm \frac{1}{2}}].
\end{align*}
So it is easy to get the contradiction if we assume there exists
$n_{m,g,Q}(\tau)$ which is not an integer.
\end{proof}

\vskip 30pt

$$ \ \ \ \ $$


\begin{thebibliography}{999}

\bibitem{AV} A. Aganagic and C. Vafa, {\em  Mirror symmetry, D-branes and
counting holomorphic discs}. arXiv: hep-th/0012041.

\bibitem{AV2}  A. Aganagic and C. Vafa,  {\em Large N Duality, Mirror Symmetry, and a Q-deformed A-polynomial
for Knots}, arXiv: 1204.4709.



\bibitem{CLPZ} Q. Chen, K. Liu, P. Peng and S. Zhu, {\em Congruent skein
relations for colored HOMFLY-PT invariants and colored Jones
polynomials}, arxiv:1402.3571v3.


\bibitem{DW} M. Dedushenko, E. Witten, {\em Some Details On The Gopakumar-Vafa
and Ooguri-Vafa Formulas}, arXiv:1411.7108.

\bibitem{EO} E. Eynard and N. Orantin, {\em Computation of open Gromov-Witten
invariants for toric Calabi- Yau 3-folds by topological recursion, a
proof of the BKMP conjecture}, Comm. Math. Phys. 337 (2015), no. 2,
483-567.

\bibitem{FGSA} H. Fuji, S. Gukov, P. Sulkowski,  H. Awata, {\em Volume Conjecture:
Refined and Categorified},  Adv. Theor. Math. Phys. {\bf 16} (2012)
1669-1777.

\bibitem{FGS} H. Fuji, S. Gukov,  P. Sulkowski, {\em Super-A-polynomial for knots
and BPS  states},  Nucl.Phys. {\bf B867} (2013) 506-546.


\bibitem{GKS} S. Garoufalidis, P. Kucharski and P. Sulkowski, {\em Knots, BPS
states, and algebraic curves}, Commun. Math. Phys. 346 (2016)
75-113.


\bibitem{GV0} R. Gopakumar, C. Vafa, {\em M-Theory and Topological
Strings-I}, arXiv:hep-th/9809187.


\bibitem{GV1} R. Gopakumar, C. Vafa, {\em M-theory and topological strings-II},
arXiv:hep-th/9812127.

\bibitem{GV2} R. Gopakumar and C. Vafa, {\em On the gauge theory/geometry
correspondence}, Adv. Theor. Math. Phys.3(5) (1999) 1415-1443.



\bibitem{HKKPTVVZ} K. Hori, S. Katz, A. Klemm, R. Pandharipande, R. Thomas, C. Vafa, R. Vakil and E. Zaslow,
{\em Mirror symmetry},(Clay mathematics monographs. 1).


\bibitem{KRSS1} P. Kucharski, M. Reineke, M. Stosic, P. Sulkowski,
{\em BPS states, knots and quivers}, arXiv:1707.02991.

\bibitem{KRSS2} P. Kucharski, M. Reineke, M. Stosic, P. Sulkowski, {\em Knots-quivers
correspondence}, arXiv:1707.04017.



\bibitem{KS0} A. Klimyk and K. Schmudgen, {\em Quantum groups and
their representation theory}, Springer-Verlag, Berlin Heidelberg
1997.


\bibitem{KS1} P. Kucharski and P. Sulkowski, {\em BPS counting for knots and combinatorics on
words}, arXiv:1608.06600.

\bibitem{LLZ} C.-C. Liu, K. Liu, J. Zhou, {\em A proof of a conjecture of Mari\~no-Vafa on Hodge integrals},
J. Differential Geom. {\bf 65}(2003).


\bibitem{LM1} J.M.F. Labastida and M. Mari\~no, {\em Polynomial invariants for torus
knots and topological strings} Comm. Math. Phys. {\bf 217}
(2001),no. 2, 423.

\bibitem{LM2} J.M.F. Labastida and M. Mari\~no, {\em A new point of view in the theory
of knot and link invariants} J. Knot Theory Ramif. {\bf 11} (2002),
173.

\bibitem{LMV} J.M.F. Labastida, M. Mari\~no and C. Vafa, {\em Knots, links and branes
at large N}, J. High Energy Phys. 2000, no. {\bf 11}, Paper 7.

\bibitem{LP} K. Liu and P. Peng, {\em Proof of the
Labastida-Mari\~no-Ooguri-Vafa conjecture}. J. Differential Geom.,
85(3):479-525, 2010.



\bibitem{LZ} W. Luo and S. Zhu, {\em Integrality structures in topological
strings I: framed unknot}, arXiv:1611.06506.


\bibitem{MMMS} A. Mironov, A. Morozov, An. Morozov, A. Sleptsov, {\em Gaussian
distribution of LMOV numbers}, arXiv:1706.00761.

\bibitem{MMMRSS} A. Mironov, A. Morozov, An. Morozov, P. Ramadevi, Vivek Kumar Singh,
A. Sleptsov, {\em Checks of integrality properties in topological
strings}, arXiv:1702.06316.




\bibitem{MV} M. Mari\~no, C. Vafa, {\em Framed knots at large N}, in:
Orbifolds Mathematics and Physics, Madison, WI, 2001, in: Contemp.
Math., vol.310, Amer. Math. Soc., Providence, RI, 2002, pp.185-204.


\bibitem{OV} H. Ooguri, C. Vafa, {\em Knot invariants and topological strings}.
Nucl. Phys. B 577(3), 419-438 (2000).



\bibitem{W} E. Witten, {\em Chern-Simons Gauge Theory As A String Theory}, Prog.
Math. {\bf 133}, 637 (1995).

\bibitem{Zhou} J. Zhou, {\em A proof of the full Mari\~no-Vafa conjecture}. Math.
Res. Lett. 17 (2010), no. 6, 1091-1099.

\bibitem{Zhu1} S. Zhu, {\em Colored HOMFLY polynomials via skein theory}, J.
High. Energy. Phys. 10(2013), 229

\bibitem{Zhu2} S. Zhu, {\em Topological strings, quiver representations and Rogers-Ramanujan
identities},  The Ramanujan Journal (2019) Vol 48(2), 399-421.

\bibitem{Zhu3} S. Zhu, {\em Topological Strings and Their
Applications in Mathematics}, Notices of The ICCM, Dec. 2017.

\end{thebibliography}
\end{document}